\documentclass[preprint,review,sort&compress,3p]{elsarticle}
\usepackage[utf8]{inputenc} 
\UseRawInputEncoding
\usepackage{mathrsfs}
\usepackage{amsfonts}
\usepackage{cancel}
\usepackage{graphicx}
\usepackage{amssymb}

\usepackage{amsthm}
\usepackage{amsmath}
\usepackage{epstopdf}
\usepackage{verbatim}
\usepackage{booktabs}%
\usepackage{flafter}
\usepackage{algorithm}
\usepackage{algorithmicx}
\usepackage{algpseudocode}
\usepackage{makecell}
\usepackage{caption,subcaption}
\usepackage{multirow}
\usepackage{appendix}
\usepackage{float}
\usepackage{array, caption, threeparttable}
\usepackage[font=small,labelfont=bf,labelsep=none]{caption}
\numberwithin{equation}{section}

\journal{ }
\begin{document}
\begin{frontmatter}
\title{Greedy double subspaces coordinate descent method via orthogonalization}
\cortext[cor1]{Corresponding author}
\author[a]{Li-Li Jin}
\author[a]{Hou-Biao Li\corref{cor1}}
\ead{lihoubiao0189@163.com}

\address[a]{School of Mathematical Sciences, University of Electronic Science and Technology of China, Chengdu, 611731, P. R. China}

\begin{abstract}
The coordinate descent method is an effective iterative method for solving large linear least-squares problems. In this paper, for the highly coherent columns case, we construct an effective coordinate descent method which iteratively projects the estimate onto a solution space formed by two greedily selected hyperplanes via Gram-Schmidt orthogonalization. Our methods may be regarded as a simple block version of coordinate descent method which involves two active columns. The convergence analysis of this method is provided and numerical simulations also confirm the effectiveness for  matrices with highly coherent columns.
\end{abstract}

\end{frontmatter}

\section{Introduction}
In practical applications, such as regression analysis and data fitting, we usually find $x\in \mathbb{R}^{n}$ to minimize the function $f(x)=\left\| b-Ax \right \|^{2} $, namely,
\begin{gather}\label{eq:1.1}
x_{*}:=\arg \min_{x\in \mathbb{R}^{n}} \left\| b-Ax \right \|^{2},
\end{gather}
where $A\in {{\mathbb{R}}^{m\times n}}$ is a full column rank matrix with $m \textgreater n $, $b\in {{\mathbb{R}}^{m}}$, and $x\in {\mathbb{R}^{n}}$ is the $n$-dimensional unknown vector. $||\cdot||$ indicates the Euclidean norm of either a vector or a matrix throughout this paper. Iterative methods such as coordinate descent method \cite{AIM,IMSL,ISLSS,CD}, which is also called the Gauss-Seidel method, are considered for solving linear least squares problem (\ref{eq:1.1}).
The coordinate descent method, starting from an initial guess $x_{0}$, can be formulated as
\begin{align}\label{eq:1.2}
x_{k+1}=x_{k}+\frac{A_{\left(j_{k}\right)}^{T}\left(b-A x_{k}\right)}{\left\|A_{\left(j_{k}\right)}\right\|^{2}} e_{j_{k}}, \quad k=0,1,2, \ldots
\end{align}
where $A_{\left(j\right)}$ is the $ j $th column of $ A $ and $(\cdot)^{T}$ represents the transpose of either a vector or a matrix. Hence $ x_{k+1} $ is obtained by projecting the current iterate $x_{k}$ to the hyperplane $H_{j_{k}}=\left\{x \mid A_{\left(j_{k}\right)}^{T} A x=A_{\left(j_{k}\right)}^{T}b \right\} $. It is  crucial to select a distinguished descent direction $e_{j_{k}}$ in the $ k $th iteration, that is, to select the column index $j_{k}$. When $j_{k}$ = ($k$ mod $n$) + 1, the coordinate descent method is the classical Gauss-Seidel method directly applied to the normal equation $ A^{T}Ax=A^{T}b$,
where $A^{T}A $ is symmetric positive definite. To speed up the convergence of the coordinate descent method, Strohmer and Vershynin \cite{RK} proposed randomized Kaczmarz method with expected exponential rate of convergence, and Leventhal and Lewis \cite{RCD} constructed the randomized coordinate descent method by randomly choosing a coordinate direction as a search direction based on an appropriate probability distribution. Soon after that, Ramdas, Needell and Ma \cite{2015REGS} provided an extended version of randomized Gauss-Seidel which converges linearly to the least norm solution in the under-determined case. In addition, Nutini et al.\cite{2015GSR} argued that in these contexts greedy selection rules gave faster convergence rates than random selection rules, and analyzed Gauss-Southwell(GS) rule and Gauss-Southwell-Lipschitz (GSL) rule. some randomized and greedy versions of the coordinate descent method can be found in \cite{CCD,GCD,RGCD,Niu}. In 2021, Liu, Jiang and Gu\cite{20212SGS} proposed a two step Gauss-Seidel (2SGS) algorithm involved two active columns by the maximum residual rule. Motivated by Needell and Ward\cite{20132SRK}, here we consider the coordinate descent method in case of the matrix $A$ with highly coherent columns, which means that the angle between these hyperplanes $H_{j}=\left\{x \mid A_{\left(j\right)}^{T} A x=A_{\left(j \right)}^{T}b \right\}(j=1,\ldots,n)$ is very small. By greedily selecting two linearly independent columns and introducing the auxiliary descent direction vector, we propose the greedy double subspaces coordinate descent method which adopts the Gram-Schmidt Orthogonalization \cite{Schmidt} process to project the current iteration onto the  corresponding solution space of two selected hyperplane.

This paper is organized as follows. In section 2, we present necessary notations. In section 3, we describe our new method and its convergence theory. In section 4, some numerical examples are provided to illustrate the effectiveness of our new methods for matrices with highly coherent columns.

\section{Notation and preliminaries}
For a vector $z\in \mathbb{R}^{n}$, $z^{(i)}$ denotes its $i\mathrm{th}$ entry. In addition, for a matrix $A$, $A_{(j)}$ represents its $j\mathrm{th}$ column. $\left\| A \right\|$, 
$\left\|A\right\|_{\infty}$ are used to indicate Euclidean norm and the infinite norm of $A$ respectively. In what follows, we denote the transpose of the matrix $A$ by $A^{T}$ and $\left\|x \right\|_{G}=\sqrt{x^{T}Gx}$ when $G$ is positive definite. $e_{j}$ is used to represent the column vector with a one in position $j$ and zero in all other positions. We use $x_{\star}=A^{\dagger}b$ to represent the unique least squares solution of (\ref{eq:1.1}). For $ x,y \in \mathbb{R}^{n} $, ${\left\langle x,y \right\rangle}=\sum_{i=1}^{n}x^{(i)}y^{(i)}$ remarks standard Euclidean inner product of two vectors.

\section{The greedy double subspaces coordinate descent method }\label{sec3}
We will assume throughout that the matrix $A$ of columns are normalized, meaning that its columns has unit Euclidean norm.
In 2015, Nutini et al.\cite{2015GSR} analyze Gauss-Southwell(GS) rule $$j_{k}=\mathop{\arg \max}\limits_{1 \leq j \leq n}\left|\nabla_{j} f\left(x_{k}\right)\right|$$
and Gauss-Southwell-Lipschitz(GSL) rule
$$j_{k}=\mathop{\arg \max}\limits_{1 \leq j \leq n} \frac{\left|\nabla_{j} f\left(x_{k}\right)\right|}{\sqrt{L_{j}}}.$$
If $f(x)=\left\| b-Ax \right \|^{2} $, and the residual vector of the normal equation $s_{k}=A^{T}(b-Ax_{k})$ is calculated in the $k$th iteration, the index $ j_{k} $ selected by GS rule satisfies
 $$j_{k}=\mathop{\arg \max}\limits_{1 \leq j \leq n}\left|s_{k}^{(j)}\right|$$
 and the index $j_{k}$ selected by GSL rule satisfies
 $$j_{k}=\mathop{\arg \max}\limits_{1 \leq j \leq n}\left \{\frac{\left|s_{k}^{(j)}\right|}{\left\|A_{(j)}\right\|}\right\}.$$
Note that GS rule is the same as GSL rule when the columns of the matrix $A$ are normalized.
This greedy coordinate descent(GCD) method with GSL rule is described by Algorithm 1.
\begin{algorithm}[H]
	\caption{The GCD Algorithm \cite{2015GSR}}
	\label{1}
	\begin{algorithmic}[1]  	
		\Require $A \in \mathbb{R}^{m \times n}, b \in \mathbb{R}^{m}, \ell,$ initial estimate $x_{0}$.
		\Ensure $x_{\ell}$.
		\For {$k=0,1,2, \ldots, \ell-1$}
		\State Compute $ s_{k}=A^{T}(b-Ax_{k}) $
		\State Select
		$j_{k}=\mathop{\arg \max}\limits_{1 \leq j \leq n}\left| s_{k}^{(j)}\right|$
		\State Update	$x_{k+1}=x_{k}+s_{k}^{\left(j_{k}\right)}e_{j_{k}}$
		\EndFor
	\end{algorithmic}
\end{algorithm}
In addition, Liu, Jiang and Gu\cite{20212SGS} construct the two step Gauss-Seidel(2SGS) algorithm by GS rule as listed Algorithm 2.
\begin{algorithm}[H]
	\caption{ The 2SGS Algorithm\cite{20212SGS}}
	\label{1}
	\begin{algorithmic}[1]  	
		\Require $A \in \mathbb{R}^{m \times n}, b \in \mathbb{R}^{m}, \ell,$ initial estimate $x_{0}$.
		\Ensure  $x_{\ell}$.
		\For {$k=0,1,2, \ldots, \ell-1$}
		\State Compute	$ s_{k}=A^{T}(b-Ax_{k}) $
		\State Select	
		$j_{k_{1}}=\mathop{\arg \max}\limits_{1 \leq j \leq n}\left| s_{k}^{(j)}\right|$, $j_{k_{2}}=\mathop{\arg \max}\limits_ {1 \leq j \leq n \atop j \neq j_{k_{1}} }\left\{\left| s_{k}^{(j)}\right|\right\}$
		\State Update	$x_{k+1}=x_{k}+s_{k}^{\left(j_{k_{1}}\right)}e_{j_{k_{1}}}+s_{k}^{\left(j_{k_{2}}\right)}e_{j_{k_{2}}}$
		\EndFor
	\end{algorithmic}
\end{algorithm}
Considering the columns of the matrix $A$ are highly coherent, which means that the angle between these hyperplanes corresponding to the normal equation is very small, the convergence speed of some versions of coordinate descent with one working column will become slow. However, if the hyperplane that is the farthest and closest to the current iteration point can be selected in $k$th iteration,
$$j_{k_{1}}=\mathop{\arg \max}\limits_{1 \leq j \leq n}\left\{\left|A_{(j)}^{T}(b-Ax_{k})\right|\right\} \quad \text { and }\quad  j_{k_{2}}=\mathop{\arg \min}\limits_{1 \leq j \leq n}\left\{\left|A_{(j)}^{T}(b-Ax_{k})\right|\right\},$$
which the angle between the two planes $ H_{j_{k_{1}}} $ and $ H_{j_{k_{2}}} $ is relatively large, it can expect to speed up the convergence of this method. Given $ x_{0} $, we first initialize $ x_{0} $ by one-step iterative GCD algorithm to get an approximate solution $ x_{1} $, while reserve the selected column index marked as $j_{0_{1}} $.
In the first iteration, since
$ A_{(j_{0_{1}})}^{T}(b-Ax_{1})=0 $, which means that the hyperplane $ H_{j_{0_{1}}} $ is the closest  to $ x_{1} $. We rewrite $ j_{0_{1}} $ as the second index of the first iteration $ j_{1_{2}} $. In addition, we can use the GSL rule to select the hyperplane farthest from the current iteration point $x_{1}$, and remark the corresponding column index as the first index of the current iteration $j_{1_{1}} $ and preserve it for next iteration. Here we construct the auxiliary descent direction vector $w_{1} $, and  $x_{2}$ can be obtained by projecting 
$x_{1}$ onto the solution space formed by these two planes 
$H_{j_{1_{1}}}$ and $H_{j_{1_{2}}}$ after some simple algebraic calculations. In the second iteration, since
$A_{(j_{1_{1}})}^{T}(b-Ax_{2})=0 $, which means that $x_{2}$ falls on the hyperplane $H_{j_{1_{1}}}$. We rewrite $j_{1_{1}}$ as the second index of the second iteration $ j_{2_{2}} $. In addition, we select the hyperplane farthest from the current iteration point $x_{2}$ by GSL rule, and denote the corresponding column index as the first index of the second iteration $j_{2_{1}}$ and keep it for next iteration.  Analogously, we introduce the auxiliary descent direction vector $w_{2} $, and get $x_{3}$ by projecting $ x_{2} $ onto the intersection of two hyperplane  $H_{j_{2_{1}}}$ and $ H_{j_{2_{2}}}$.
Sequentially, we gain the sequence $ \left\{x_{k}\right\}_{k=0}^{\infty} $in the same way. The specific algorithm is detailedly exhibited in the Algorithm 3.

\begin{algorithm}[H]
	\caption{The greedy double subspaces coordinate descent(GDSCD) method}
	\label{2}
	\begin{algorithmic}[1]  
		\Require $A \in \mathbb{R}^{m \times n}, b \in \mathbb{R}^{m},\; \ell\; $, initial estimate $x_{0}$.
		\Ensure  $x_{\ell}$.
		\State Compute $ s_{0}=A^{T}(b-Ax_{0}),
		j_{0_{1}}=\mathop{\arg \max}\limits_{1 \leq j \leq n}\left\{\left| s_{0}^{(j)}\right| \right\}$
		\State Update $ x_{1}=x_{0}+s_{0}^{(j_{0_{1}})}e_{j_{0_{1}}} $
		\For {$k=1,2, \ldots, \ell-1$}
		\State Select
		$j_{k_{1}}=\mathop{\arg \max}\limits_{1 \leq j \leq n}\left\{\left| s_{k}^{(j)}\right| \right\}$, $ j_{k_{2}}=\mathop{\arg \max}\limits_ {1 \leq j \leq n}\left\{\left| s_{k-1}^{(j)}\right|\right\}$
		\State  Set
		$\widetilde{{\mu}}_{k}=\left\langle A_{(j_{k_{1}})},A_{(j_{k_{2}})}\right\rangle$,
		$\widetilde{y}_{k}=x_{k}+s_{k}^{(j_{k_{1}})}e_{j_{k_{1}}}$
		\State  Set
		$\widetilde{v}_{k}=\dfrac{A_{(j_{k_{2}})}-{\widetilde{\mu}_{k}}A_{(j_{k_{1}})}}{\sqrt{1-{\widetilde{\mu}_{k}}^{2}}}$,
		$\widetilde{\beta}_{k}=\dfrac{(A_{(j_{k_{2}})}-{\widetilde{\mu}_{k}}A_{(j_{k_{1}})})^{T}b}{\sqrt{1-{\widetilde{\mu}_{k}}^{2}}}$,
		$w_{k}=\dfrac{e_{j_{k_{2}}}-{\widetilde{\mu}_{k}}e_{j_{k_{1}}}}{\sqrt{1-{\widetilde{\mu}_{k}}^{2}}} $
		\State  Update	$x_{k+1}=\widetilde{y_{k}}+  \left(\widetilde{\beta}_{k}-\widetilde{v}_{k}^{T}A\widetilde{y}_{k}\right)w_{k}$
		\EndFor
	\end{algorithmic}
\end{algorithm}
For the convergence theory of the GDSCD algorithm, we have the following theorem.

\noindent\textbf{Theorem 1.} The iteration sequence $\left\{x_{k}\right\}_{k=0}^{\infty}$, generated by the GDSCD method starting from any initial guess $x_{0}\in {\mathbb{R}}^{n} $, linearly converges to the unique least-squares solution $x_{\star}=A^{\dagger}b$ and satisfies
$$ \left\|x_{*}-x_{1}\right\|_{A^{T}A}^{2}\leq \left(1-\sigma(A^{T},\infty)^{2}\right)\left\|x_{*}-x_{0}\right\|_{A^{T}A}^{2},$$
and
$$\left\|x_{*}-x_{k+1}\right\|_{A^{T}A}^{2} \leq \left(1-\dfrac{1}{1-\widetilde{{\mu}}_{k}^{2}} \sigma(A^{T},\infty)^{2}\right)\left\|x_{*}-x_{k}\right\|_{A^{T}A}^{2},\quad k=1,2,\dots ,$$
where
$ \widetilde{\mu}_{k}=\left\langle A_{(j_{k_{1}})},A_{(j_{k_{2}})}\right\rangle$, and $ \sigma(A^{T},\infty)$ is  the Hoffman-like constant as defined in \cite{2016Orthogo}.
\vspace{1mm}

\noindent\textbf{Proof.}
Denote $ r_{k}=b-Ax_{k}, s_{k}=A^{T}r_{k}, e_{k}=x^{*}-x_{k}$.
From Algorithm 3, for $ k\geq 1 $, we have
\begin{align*}
x_{k+1} &=\widetilde{y}_{k}+(\widetilde{\beta}_{k}-\widetilde{v}_{k}^{T}A\widetilde{y}_{k})w_{k}\\
&=x_{k}+s_{k}^{(j_{k_{1}})}e_{j_{k_{1}}} +(\widetilde{\beta}_{k}-\widetilde{v}_{k}^{T}A\widetilde{y}_{k})w_{k}.
\end{align*}
Hence, we can get 
$$A(x^{*}-x_{k+1})=A(x^{*}-x_{k})-s_{k}^{(j_{k_{1}})}A_{(j_{k_{1}})}-
(\widetilde{\beta}_{k}-\widetilde{v}_{k}A\widetilde{y}_{k})\widetilde{v}_{k}.$$
Since $\widetilde{v}_{k} $ is
perpendicular to $A_{(j_{k_{1}})}$ for any $ k\geq 1 $, it holds that
\begin{align}
\left\|Ae_{k+1}\right\|^{2}
&=\left\|Ae_{k}\right\|^{2}-\left(s_{k}^{(j_{k_{1}})}\right)^{2}+(\widetilde{\beta}_{k}-\widetilde{v}_{k}^{T}A\widetilde{y}_{k})^{2}-2(\widetilde{\beta}_{k}-\widetilde{v}_{k}^{T}A\widetilde{y}_{k})\left\langle \widetilde{v}_{k},Ae_{k}\right \rangle.\label{(3.1)}
\end{align}
we can easily calculate
\begin{align}
\widetilde{\beta}_{k}-\widetilde{v}_{k}^{T}A\widetilde{y}_{k}&=\widetilde{\beta}_{k}-\widetilde{v}_{k}^{T}A\widetilde{x}_{k} \notag \\
&=\dfrac{A_{(j_{k_{2}})}^{T}b-A_{(j_{k_{2}})}^{T}Ax_{k}-\widetilde{\mu}_{k}(A_{(j_{k_{1}})}^{T}b-A_{(j_{k_{1}})}^{T}Ax_{k})}{\sqrt{1-\widetilde{\mu}_{k}^{2}}} \notag \\
&=\dfrac{A_{(j_{k_{2}})}^{T}-\widetilde{\mu}_{k}A_{(j_{k_{1}})}^{T}}{\sqrt{1-\widetilde{\mu}_{k}^{2}}}A(x^{*}-x_{k}) \notag \\
&=\widetilde{v}_{k}^{T}Ae_{k}.\label{(3.2)}
\end{align}
For $ k\geq 2 $
\begin{align*}
A_{(j_{k-1_{1}})}^{T}Ax_{k}&=A_{(j_{k-1_{1}})}^{T}A\left(x_{k-1}+s_{k-1}^{(j_{k-1_{1}})}e_{j_{k-1_{1}}} +(\widetilde{\beta}_{k-1}-\widetilde{v}_{k-1}^{T}A\widetilde{y}_{k-1})w_{k-1} \right) \notag \\
&=A_{(j_{k-1_{1}})}^{T}Ax_{k-1}+s_{k-1}^{(j_{k-1_{1}})}A_{(j_{k-1_{1}})}^{T}A_{(j_{k-1_{1}})}+(\widetilde{\beta}_{k-1}-\widetilde{v}_{k-1}^{T}A\widetilde{y}_{k-1})A_{(j_{k-1_{1}})}^{T} \widetilde{v}_{k-1} \notag \\
&=A_{(j_{k-1_{1}})}^{T}Ax_{k-1}+s_{k-1}^{(j_{k-1_{1}})}  \\
&=A_{(j_{k-1_{1}})}^{T}b,
\end{align*}
for $ k=1 $,
	\begin{align*}
A_{(j_{0_{1}})}^{T}Ax_{1}&=A_{(j_{0_{1}})}^{T}A(x_{0}+s_{0}^{(j_{0_{1}})}e_{j_{0_{1}}}) \\
&=A_{(j_{0_{1}})}^{T}Ax_{0}+s_{0}^{(j_{0_{1}})}A_{(j_{0_{1}})}^{T}A_{(j_{0_{1}})}\\
&=A_{(j_{0_{1}})}^{T}Ax_{0}+A_{(j_{0_{1}})}^{T}(b-Ax_{0}) \\
&=A_{j_{0_{1}}}^{T}b,
\end{align*}
which implies that $ s_{k}^{(j_{k-1_{1}})}=0$ for any $k\geq 1$. Therefore, we can simplify this equation (\ref{(3.2)})
\begin{align}
\widetilde{v_{k}}^{T}Ae_{k} &=\dfrac{A_{(j_{k_{2}})}^{T}-\widetilde{\mu}_{k}A_{(j_{k_{1}})}^{T}}{\sqrt{1-\widetilde{{\mu}}_{k}^{2}}}Ae_{k} \notag \\
&=\dfrac{s_{k}^{(j_{k_{2}})}-\widetilde{\mu}_{k}s_{k}^{(j_{k_{1}})}}{\sqrt{1-\widetilde{\mu}_{k}^{2}}}Ae_{k} \notag \\
&=\dfrac{s_{k}^{(j_{k-1_{1}})}-\widetilde{\mu}_{k}s_{k}^{(j_{k_{1}})}}{\sqrt{1-\widetilde{\mu}_{k}^{2}}}Ae_{k} \notag \\
&=\dfrac{-\widetilde{\mu}_{k}s_{k}^{(j_{k_{1}})}}{\sqrt{1-\widetilde{\mu}_{k}^{2}}}Ae_{k}. \label{(3.3)}
\end{align}	
Thus, combining (\ref{(3.1)}) and (\ref{(3.3)}), we can further obtain
\begin{align}
\left\|Ae_{k+1}\right\|^{2}
&=\left\|Ae_{k}\right\|^{2}-\left(s_{k}^{(j_{k_{1}})}\right)^{2}+(\widetilde{\beta}_{k}-\widetilde{v}_{k}^{T}A\widetilde{y}_{k})^{2}-2(\widetilde{\beta}_{k}-\widetilde{v}_{k}^{T}A\widetilde{y}_{k})\left\langle \widetilde{v}_{k},Ae_{k}\right \rangle \notag \\
&=\left\|Ae_{k}\right\|^{2}-{s_{k}^{(j_{k_{1}})}}^{2}-\dfrac{\widetilde{{\mu}}_{k}^{2}}{1-\widetilde{\mu}_{k}^{2}}{s_{k}^{(j_{k_{1}})}}^{2} \notag \\
&=\left\|Ae_{k}\right\|^{2}-\dfrac{1}{1-\widetilde{\mu}_{k}^{2}}{s_{k}^{(j_{k_{1}})}}^{2} \notag \\
&=\left\|Ae_{k}\right\|^{2}-\dfrac{1}{1-\widetilde{\mu}_{k}^{2}} \mathop{ \max}\limits_{i\in[n]}\left\{ \left| s_{k}^{(j_{k_{1}})^{2}} \right|\right\}  \notag \\
&=\left\|Ae_{k}\right\|^{2}-\dfrac{1}{1-\widetilde{\mu}_{k}^{2}} \left\|A^{T} A(x^{*}-x_{k})\right\|_{\infty}^{2}  \notag \\
&\leq \left\|Ae_{k}\right\|^{2}-\dfrac{1}{1-\widetilde{\mu}_{k}^{2}} \sigma(A^{T},\infty)^{2}\left\|Ae_{k}\right\|^{2}  \label{(3.4)}\\
&=(1-\dfrac{1}{1-\widetilde{\mu}_{k}^{2}} \sigma(A^{T},\infty)^{2})\left\|Ae_{k}\right\|^{2}, \label{(3.5)}
\end{align}
where the equation (\ref{(3.4)}) is achieved with the use of Hoffman-like constant $ \sigma(A^{T},\infty)$ in \cite{2016Orthogo}.

In addition, for $ k=0 $, we have
\begin{align}
\left\|A(x_{*}-x_{1})\right\|^{2}&=\left\|A(x_{*}-x_{0})\right\|^{2}-\left(s_{0}^{(j_{0_{1}})}\right)^{2} \notag\\
&=\left\|A(x_{*}-x_{0})\right\|^{2}-\left\|A^{T}A(x_{*}-x_{0})\right\|_{\infty}^{2} \notag \\
&=\left\|A(x_{*}-x_{0})\right\|^{2}-\dfrac{\left\|A^{T}A(x_{*}-x_{0})\right\|_{\infty}^{2}}{\left\|A(x_{*}-x_{0})\right\|^{2}}\left\|A(x_{*}-x_{0})\right\|^{2} \notag \\
&\leq (1-\sigma(A^{T},\infty)^{2})\left\|A(x_{*}-x_{0})\right\|^{2},\label{(3.6)}
\end{align}
where the equation (\ref{(3.6)}) is obtained by making use of the Hoffman-like constant $\sigma(A^{T},\infty)$ in \cite{2016Orthogo}. Therefore, we can get the estimate by combining the equation (\ref{(3.5)}) and (\ref{(3.6)}).

\section{Numerical experiments}
\vspace{-2mm}
In this section we perform several experiments to compare the convergence rate of the GCD, 2SGS and GDSCD methods for some column normalized matrices in terms of the number of iteration steps (denoted as "IT") and the computing time in seconds (denoted as "CPU"). Note that the IT and CPU listed in our numerical results denote the arithmetical averages of the required iteration numbers and the CPU times with respect to 30 times repeated runs of the corresponding methods. To be precise, we define the pairwise coherence of a column standardized matrix as follows by reference to the notations in \cite{20132SRK}
$$\delta=\delta(A)=\min _{i \neq j}|\langle A_{(i)}, A_{(j)}\rangle| \quad \text{and} \quad \Delta=\Delta(A)=\max _{i \neq j}|\langle A_{(i)}, A_{(j)}\rangle|.$$
In addition, we construct various types of 500 $ \times $ 100 and 5000 $ \times $ 500 matrices $A$ by MATLAB function $\mathit{unifrnd}$ which generate continuous uniformly distributed random variables on the interval $[c,1]$. Thus matrices with different pairwise coherence  can be constructed by changing the value of $c $. Moreover, the corresponding linear system  \eqref{eq:1.1} could be either consistent or inconsistent. The solution vector $x_{*}$ is generated by using MATLAB function $\mathit{randn}$. As for $b$, we take $b=A x_{*}$ when the systems is consistent, and $b=Ax_{*}+b_{0}$ when the linear system is inconsistent, where $b_{0}$ is a nonzero vector belonging to the null space of $A^{T}$, and null($A^{T}$) is generated by using MATLAB function $\mathit{null}$.
\begin{table}[H]
	\centering
	\setlength{\abovecaptionskip}{5pt}%
	\setlength{\belowcaptionskip}{10pt}%
	\caption{\;\; IT and CPU for various $ 500\times 100 $ matrices when the linear system is consistent.}
\begin{tabular}{cccccccc}
	\toprule
	$ 500\times100$ & & $A_{1}$ & $A_{2}$ & $A_{3}$ & $A_{4}$ & $A_{5}$& $A_{6}$ \\
	\hline
	$ c $ & \ &-0.8 & -0.1 & 0.8 & 0.85& 0.9& 0.95 \\
	$ \delta $&\ &4.0558e-06&0.6059&0.9951&0.9975&0.9989&0.9997 \\
	$ \Delta $&\ &0.1890&0.7318&0.9966&0.9982&0.9993 &0.9998\\
	rank(A)&\ &100&100&100&100&100&100 \\
	GCD &IT&$ 494$&$ 1311 $&$ 92067 $&$ 197026 $&$ -- $ &$ -- $\\
	&CPU&$0.0270$&$ 0.0723 $&$ 5.7447 $&$ 10.8633 $&$ -- $&$ -- $ \\
	2SGS&IT&252&237&2262&4243&8768&40647\\
	&CPU&0.0170&0.0151&0.1357&0.2699&0.5342&2.3530\\
	GDSCD &IT&433&365&383&385&377&389\\
	&CPU&0.0393&0.0311&0.0332&0.0300&0.0309&0.0333\\
	\bottomrule
\end{tabular}
	\label{tablea}
\end{table}
\begin{table}[H]
	\centering
	\setlength{\abovecaptionskip}{5pt}%
	\setlength{\belowcaptionskip}{10pt}%
	\caption{\;\;IT and CPU for various $ 5000\times 500 $ matrices when the linear system is consistent.}
	\begin{tabular}{cccccccc}
		\toprule
		$ 5000\times500$ & & $B_{1}$ & $B_{2}$ & $B_{3}$ & $B_{4}$ & $B_{5}$& $B_{6}$ \\
		\hline
		$ c $ & \ &-0.8 & -0.1 & 0.8 & 0.85& 0.9& 0.95 \\
		$ \delta $&\ &2.6205e-06&0.6385&0.9956&0.9976&0.9991&0.9998 \\
		$ \Delta $&\ &0.0946&0.6943&0.9962&0.9980&0.9990 &0.9998\\
		rank(A)&\ &500&500&500&500&500&500 \\
		GCD &IT&$ 1755$&$ 5001 $&$ -- $&$ -- $&$ -- $ &$ -- $\\
		&CPU&$2.6047$&$ 8.3252 $&$ -- $&$ -- $&$ -- $&$ -- $ \\
		2SGS&IT&803&800&2682&3758&11059&39232\\
		&CPU&1.2081&1.4180&4.0388&5.6522&16.6704&40.4757\\
		GDSCD&IT&1553&1731&2033&2043&2040&2050\\
		&CPU&2.5811&3.2864&3.4455&3.2922&3.4014&3.3556\\
		\bottomrule
	\end{tabular}
	\label{tablea}
\end{table}
In all implementations, the initial point $x_{0}$ is set to be a zero vector and terminated once the relative solution error (RSE), defined by
$\text{RSE}=\dfrac{\left\| x_{k}-x_{*} \right\|^{2}}{\left\| x_{*} \right\|^{2}}$ at the current iterate $x_{k}$, satisfies RSE $\leq  10^{-6} $, or the number of iteration steps exceeds $200,000$. The latter is given a label $ ``--" $in the numerical tables. All experiments are performed on a PC with Intel(R) Core(TM) i5-8250U CPU @ 1.60GHz 8.00GB using MATLAB R2016a.

For the consistent systems, Tables 1-2 show the number of iteration steps and the computing time of various $ 500\times 100 $ and $ 5000\times 500 $ matrices for the three methods. Here the matrix $ A_{5}, A_{6},B_{5},B_{6}$ have highly coherent columns with $\delta \approx \Delta$ or $\delta =\Delta$ . Note that when the value of $c$ is closer to 1, both $\delta$ and $\Delta$ are closer to 1, which means that the columns of these matrices are more coherent. Admittedly, the
\begin{figure}[htbp]
	\centering
	\includegraphics[scale=0.65]{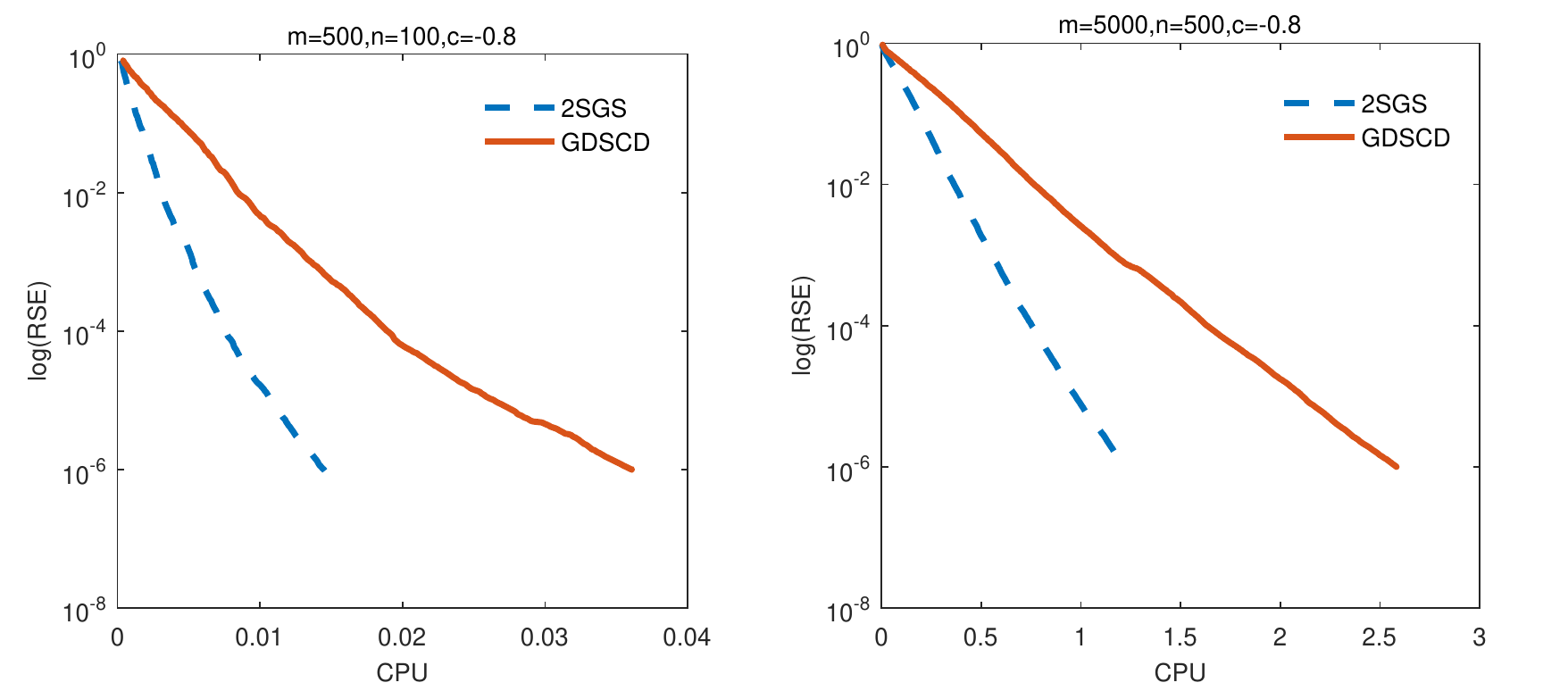}
	\caption{\;$\log_{10}(\mathrm{RSE})$ versus CPU for $ 500\times100$ (left), $ 5000\times500$ (right)  matrices when $ c $ is -0.8.}
\end{figure}
\begin{figure}[htbp]
	\centering
	\includegraphics[scale=0.7]{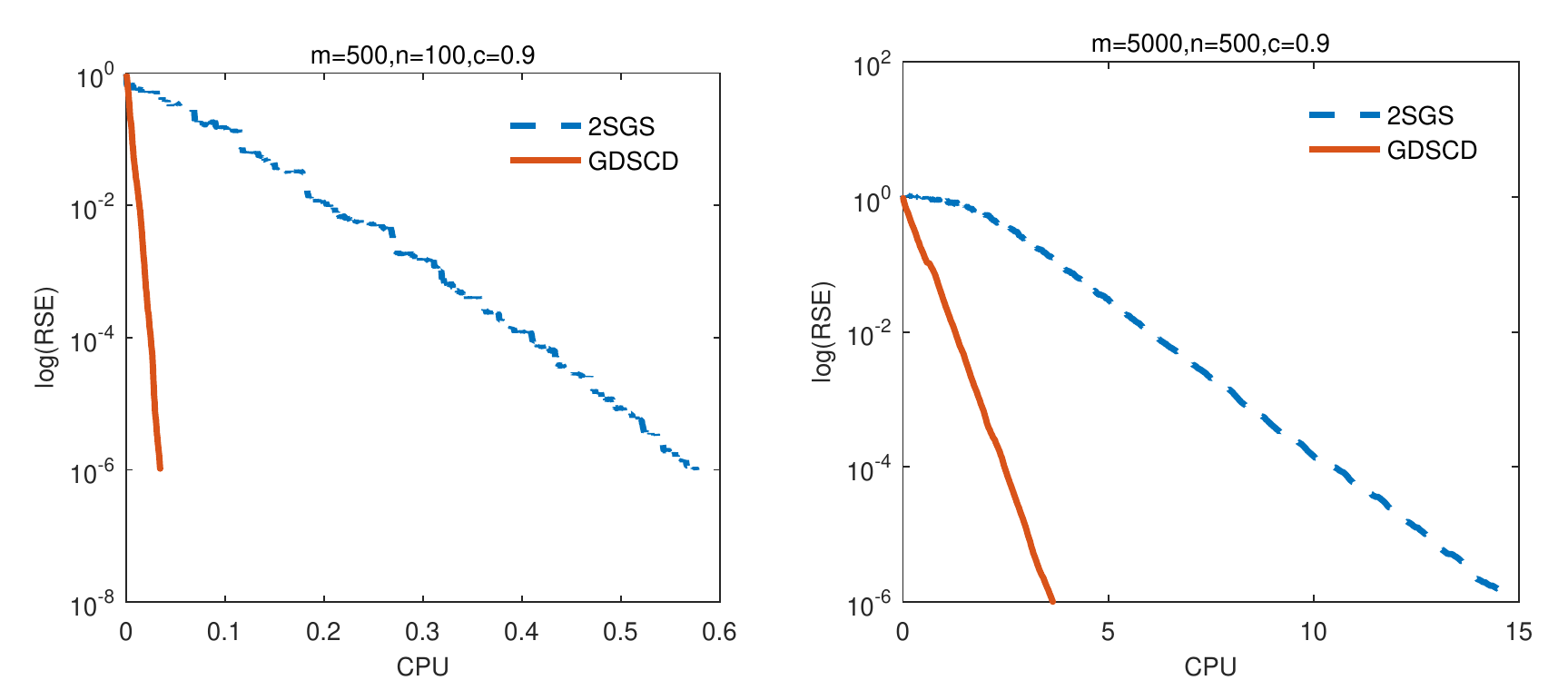}
	\caption{\;$\log_{10}(\mathrm{RSE})$ versus CPU for $ 500\times100$ (left), $ 5000\times500$ (right) matrices when $ c $ is 0.9.}
\end{figure}
convergence rate of the 2SGS method is faster than both the GCD method and the GDSCD method in the case where $ \delta$ and $\Delta $ are small(see Table1 of $ A_{1}$ and $A_{2}$, Table2 of $B_{1}$ and 
$B_{2} $). However, it is true that GDSCD method outperforms the other two algorithms in both time and number of iteration steps with increasing the coherence of column of matrices(see Table1 of $ A_{5}$ and $ A_{6}$, Table2 of $B_{5}$ and $B_{6} $). What is better is that the GDSCD method is at least 5 times
faster than the other two methods for matrices highly coherent columns(see Table1 of $ A_{6}$, Table2 of $B_{6} $). In addition, Figures 1 and 2  depict the
curves of the $\log_{10}(\mathrm{RSE})$ versus CPU time of when the linear system is
consistent with $ c=-0.8  $ and $ c=0.9 $, which further shows the effectiveness of the proposed new method for highly column-coherent matrices.
\begin{table}[H]
	\centering
	\setlength{\abovecaptionskip}{5pt}%
	\setlength{\belowcaptionskip}{10pt}%
	\caption{\;\;IT and CPU for various $ 500\times 100 $ matrices when the linear system is inconsistent.}
	\begin{tabular}{cccccccc}
		\toprule
		$ 500\times100$ & & $A_{1}$ & $A_{2}$ & $A_{3}$ & $A_{4}$ & $A_{5}$& $A_{6}$ \\
		\hline
		$ c $  & \ &-0.8 & -0.1 & 0.8 & 0.85& 0.9& 0.95 \\
		$ \delta $&\ &4.0558e-06&0.6059&0.9951&0.9975&0.9989&0.9997 \\
		$ \Delta $&\ &0.1890&0.7318&0.9966&0.9982&0.9993 &0.9998\\
		rank(A)&\ &100&100&100&100&100&100 \\
		GCD &IT&$ 465$&$ 1261 $&$ 105668 $&$ 188781 $&$ -- $ &$ -- $\\
		&CPU&$0.0281$&$ 0.0705 $&$ 5.8241 $&$ 10.3324 $&$ -- $&$ -- $ \\
		2SGS&IT&243&219&2279&3518&8378&38576\\
		&CPU&0.0167&0.0147&0.1469&0.2153&0.5080&2.2868\\
		GDSCD&IT&443&360&410&384&400&388\\
		&CPU&0.0338&0.0291&0.0329&0.0332&0.0327&0.0322\\
		\bottomrule
	\end{tabular}
	\label{tablea}
\end{table}
Similarly, Tables 3 and Tables 4 demonstrate the same conclusion for the inconsistent systems. We observe that when the column coherence of the matrix is not particularly large, the convergence efficiency of the 2SGS  method is better than the other two methods (see Table3 of $ A_{1}$, $ A_{2}$, Table4 of $B_{1}$, $ B_{2}$ ).
\begin{table}[H]
	\centering
	\setlength{\abovecaptionskip}{5pt}%
	\setlength{\belowcaptionskip}{10pt}%
	\caption{\;\;IT and CPU for various $ 500\times 100 $ matrices when the linear system is inconsistent.}
	\begin{tabular}{cccccccc}
		\toprule
		$ 5000\times500$ & & $B_{1}$ & $B_{2}$ & $B_{3}$ & $B_{4}$ & $B_{5}$& $B_{6}$ \\
		\hline
		c & \ &-0.8 & -0.1 & 0.8 & 0.85& 0.9& 0.95 \\
		$\delta$&\ &2.6205e-06&0.6385&0.9956&0.9976&0.9991&0.9998 \\
		$ \Delta$&\ &0.0946&0.6943&0.9962&0.9980&0.9990 &0.9998\\
		rank(A)&\ &500&500&500&500&500&500 \\
		GCD &IT&$ 1430$&$ 4876 $&$ -- $&$ -- $&$ -- $ &$ -- $\\
		&CPU&$2.5931$&$ 8.2934 $&$ -- $&$ -- $&$ -- $&$ -- $ \\
		2SGS&IT&789&787&2676&3821&11068&39247\\
		&CPU&1.1994&1.3290&4.0472&5.8231&16.6803&41.0543\\
		GDSCD&IT&1550&1723&2042&2109&2036&2061\\
		&CPU&2.4992&3.1263&3.4337&3.1879&3.4029&3.3421\\
		\bottomrule
	\end{tabular}
	\label{tablea}
\end{table}
\begin{figure}[htbp]
	\centering
	\includegraphics[scale=0.65]{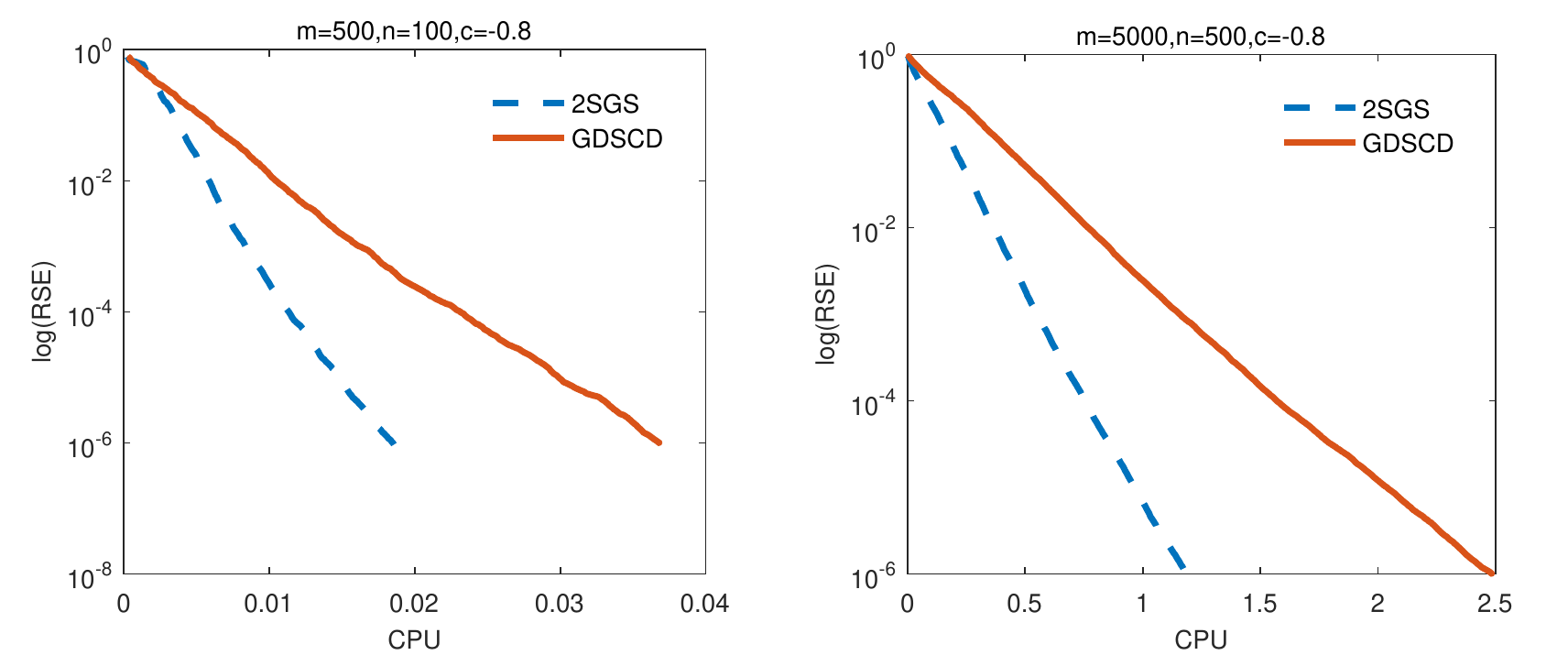}
	\caption{\;$\log_{10}(\mathrm{RSE})$ versus CPU for $500\times100$ (left), $5000\times500$ (right) matrices when $c$ is -0.8.}
\end{figure}
\begin{figure}[H]
	\centering
	\includegraphics[scale=0.65]{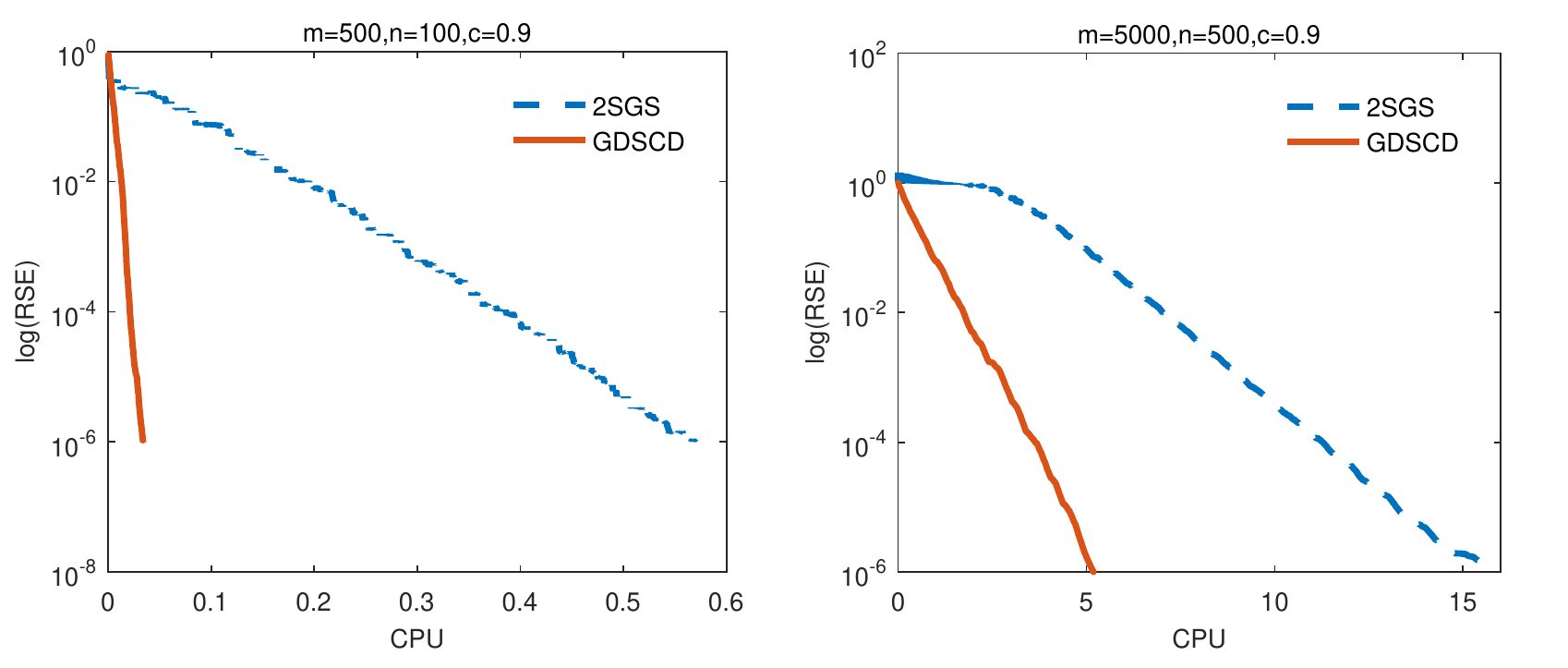}
	\caption{\;$\log_{10}(\mathrm{RSE})$ versus CPU for $ 500\times100$ (left), $5000\times500$ (right) matrices when $c$ is 0.9.}
\end{figure}
However, as the value of c getting closer to 1, that is, the coherence of columns of the constructed matrix increasing continuously, we see again that the GDSCD method significantly outperforms the GCD
method and the 2SGS method for both iteration counts and CPU times(see Table3 of $ A_{5}$, $ A_{6}$, Table4 of $B_{5}$, $ B_{6}$ ). In addition, we plot the
curves of $\log_{10}(\mathrm{RSE})$ versus CPU time of for $ 500\times 100$ and $5000\times 500$ matrices with $c=-0.8$ and $c=0.9$ when the linear system inconsistent.
\section{Conclusions}
This paper introduces a greedy block version of the coordinate descent method for the linear least-squares problem. The main idea of our method is to greedily select two linearly independent columns at first, and then projects the estimate onto a solution space of the two selected hyperplanes by Gram-Schmidt orthogonalization process. The theoretical analysis and numerical results show that our proposed method is efficient for matrices with highly coherent columns.

\textbf{Acknowledgements}. \emph{The authors are very much indebted to the referees for their constructive comments and valuable suggestions. The authors are partially supported by National Natural Science Foundation of China (11101071, 11271001, 51175443)}.
\vspace{-6mm}
\section*{}

\end{document}